\documentclass[11pt]{article}
\usepackage[colorlinks,
            linkcolor=blue,
            anchorcolor=green,
            citecolor=red
            ]{hyperref}
\usepackage{mathrsfs,amssymb,amsfonts}
\usepackage{amsmath,amssymb,latexsym,color}
\usepackage[mathscr]{eucal}
\newtheorem{thm}{Theorem}[section]
\newtheorem{defin}[thm]{Definition}
\newtheorem{prop}[thm]{Proposition}
\newtheorem{lemma}[thm]{Lemma}
\newtheorem{cor}[thm]{Corollary}
\newtheorem{example}{Example}
\newtheorem{obser}{Observation}
\newenvironment{definition}{\begin{defin} \rm}{\end{defin}}
\newtheorem{remark}{Remark}[section]

\newenvironment{eg}{\begin{example} \rm}{\end{example}}
\newenvironment{ob}{\begin{obser} \rm}{\end{obser}}

\newcommand{\proof}{{\it Proof.\quad}}
\newcommand{\qed}{\hfill\Box\medskip}
\usepackage{CJK}
\begin{document}
\begin{CJK*}{GBK}{song}
\renewcommand{\abovewithdelims}[2]{
\genfrac{[}{]}{0pt}{}{#1}{#2}}
%%%%%%%%%%%%%%%%%%%%%%%%%%%%%%%%%%%%%%%%%%%%%%%%%%%%%%%%%%%%%%%%%%%%%%%%%%%%%%%%%%%%%%%%
%%%%%%%%%%%%%%%%%%%%%%%%%%%%%%%%%%%%%%%%%%%%%%%%%%%%%%%%%%%%%%%%%%%%%%%%%%%%%%%%%%%%%%%%

\title{\bf On the power graph of a finite group}

\author{Min Feng \quad Xuanlong Ma \quad Kaishun Wang\footnote{ Corresponding author.\newline
{\em E-mail address:}  fgmn\_1998@163.com (M. Feng), mxl881112@126.com (X. ma), wangks@bnu.edu.cn (K. Wang).}\\
{\footnotesize   \em  Sch. Math. Sci. {\rm \&} Lab. Math. Com. Sys.,
Beijing Normal University, Beijing, 100875,  China} }
 \date{}
 \maketitle

\begin{abstract}
The power graph $\mathcal P_G$ of a finit group $G$ is the graph with the vertex set   $G$, where two elements are adjacent if one is a power of the other.  We first show that $\mathcal P_G$ has an transitive orientation,  so it is a perfect graph and its core is a complete graph. Then we use the poset on all cyclic subgroups (under usual inclusion) to characterise the structure of $\mathcal P_G$. Finally, the closed formula for the metric dimension of $\mathcal P_G$ is established. As an application, we compute the metric dimension of the power graph of a cyclic group.

\medskip


\end{abstract}

\section{Introduction}
In this paper, a graph means an undirected simple graph and
a digraph means a directed graph without loops.
We always use $V(\Gamma)$ and $E(\Gamma)$ to denote the vertex set
and the edge set (resp. the arc set) of a graph (resp. digraph)
$\Gamma$, respectively.
All groups, graphs and digraphs considered are finite.

Given a group,
there are different ways to associate a directed or undirected graph
to the group: intersection graphs  \cite{bos,zel},
 commuting graphs \cite{bat},  prime graphs \cite{liy}
and of course Cayley (di)graphs, which have a long history.

Let $G$ be a group. The {\em power digraph}
 of $G$ is the digraph
$\overrightarrow{\mathcal P}_G$ with the vertex set $G$, where
there is an arc from $x$ to $y$ if $x\neq y$ and $y = x^m$
for some positive integer $m$.
The {\em  power graph} $\mathcal P_G$ have the vertex
set $G$ and two distinct elements $x$ and $y$ are adjacent if one
is a power of the other.
The power digraph was introduced
by Kelarev and Quinn \cite{kel1,kel2}  and they
called it directed power graph and defined it on semigroups.
Motivated by this,
Chakrabarty, Ghosh and Sen \cite{cha} defined
 power graphs of semigroups.
Recently, Many interesting results on the power graphs have been obtained, see \cite{cam,came,doo,mir,mog,tam}.
In \cite{aba},  Abawajy,  Kelarev and Chowdhury give a survey of the current state of knowledge on this research direction by presenting all results and open questions recorded in the literature dealing with power graphs.

Given a graph $\Gamma$,
the digraph $\mathcal{O}$ is an {\em orientation}
for $\Gamma$ if $V(\mathcal{O})=V(\Gamma)$ and
$|\{(u,v),(v,u)\}\cap E(\mathcal{O})|=1$ for all $\{u,v\}\in E(\Gamma)$.
A {\em transitive orientation} for $\Gamma$ is an orientation $\mathcal{O}$
 such that $\{(u,v),(v,w)\}\subseteq E(\mathcal{O})$
implies $(u,w)\in E(\mathcal{O})$.
A {\em comparability graph} is a graph that admits a
transitive orientation. It has been originally studied in \cite{gho} and
characterized in  \cite{gal,gil}.
Recently, comparability graphs have been used to model
optimization problems in railways: see \cite{di} for a survey.
Comparability graphs have an important role in graph theory because of their
relationship with partially ordered sets:
a comparability graph is a  graph
which has the vertex set a poset and join two distinct
elements if they
are comparable in the poset.

For a graph $\Gamma$, let $d_{\Gamma}(u,v)$
denote the distance between two vertices $u$ and  $v$.
By an ordered set
of vertices, we mean a set $W=\{w_{1},\ldots,w_{k}\}$ on which the
ordering $(w_{1},\ldots,w_{k})$ has been imposed.
For an ordered
subset $W=\{w_{1},\ldots,w_{k}\}$, write
$\mathcal{D}_{\Gamma}(v|W)=(d_{\Gamma}(v,w_{1}),\ldots,d_{\Gamma}(v,w_{k}))$.
A \emph{resolving set} of $\Gamma$ is an ordered subset of
vertices $W$ such that $\mathcal{D}_{\Gamma}(u|W)=\mathcal{D}_{\Gamma}(v|W)$ if and
only if $u=v$.
 The \emph{metric dimension} of $\Gamma$, denoted by $\dim(\Gamma)$, is
the minimum cardinality of a resolving set of $\Gamma$. Metric
dimension was first introduced in the 1970s, independently by Harary
and Melter \cite{Ha} and by Slater \cite{Sl}. It is a parameter that
has appeared in various applications (see \cite{RP,Ca} for more
information).  It was noted in \cite[p. 204]{Ga} and \cite{Kh} that
 determining the metric dimension of a graph is an NP-complete problem.

In this paper, we study the power graph of a group $G$.
In Section $2$, we first construct a transitive orientation for
$\mathcal P_G$,  then get some properties of $\mathcal P_G$,
 and finally characterise the structure of
 $\mathcal P_G$ by using the poset on all cyclic subgroups of $G$ (under usual inclusion).
In Section $3$, we establish a closed formula for the metric dimension of  $\mathcal P_G$.

\section{Properties and characterization}
In this section, we get some properties of the power graph
 of a group $G$ and characterize the structure of $\mathcal P_G$.
In Subsection $2.1$, we construct a transitive orientation for  $\mathcal P_G$,
which is a subdigraph of the power digraph $\overrightarrow{\mathcal P}_G$.
Therefore, we know that
$P_G$ is a comparability graph. Then we show that it is a perfect graph
and its core is complete. Since a transitive orientation uniquely
determine a partially ordered set (or poset for simplify),
Subsection $2.2$ reviews some definitions or properties associated
with posets.
In Subsection $2.3$,
we characterize the structure of $\mathcal P_G$ by using the poset
on all cyclic subgroups of $G$ (under usual inclusion).

\subsection{Transitive orientations and comparability graphs}
Let $G$ be a group.
For $x\in G$, denote  by $[x]$ the set of all generators of the cyclic subgroup $\langle x\rangle$.
Write
\begin{equation}\label{1}
\mathcal C'(G)=\{[x]\mid x\in G\}=\{[x_{1,1}], \ldots,[x_{k,1}]\},\textup{ where } [x_{i,1}]=\{x_{i,1},\ldots,x_{i,s_i}\}.
\end{equation}
We impose an ordering $(x_{i,1},\ldots,x_{i,s_i})$ on the set $[x_{i,1}]$ for each $i\in\{1,\ldots,k\}$.
\begin{definition}\label{defin}
 For elements $x$ and $y$ in a group $G$, define $x\prec y$
if one of the followings holds.

{\rm(i)} For some $i$, $x=x_{i,l}$, $y=x_{i,t}$ and $l<t$.

{\rm(ii)} $\langle x\rangle\subsetneqq\langle y\rangle$.\\
Define $x\preceq y$ if $x\prec y$ or $x=y$.
\end{definition}

The proof of  the following lemma is clear from the above definition.
\begin{lemma}\label{lemma}
Suppose $G$ is a group.
With reference to {\rm(\ref{1})},
if there exist two distinct indices $i$ and $j$ such that
$x_{i,l_0}\prec x_{j,t_0}$ for some positive integers
$l_0$ and $t_0$, then $x_{i,l}\prec x_{j,t}$ for
each $x_{i,l}\in[x_{i,1}]$ and each $x_{j,t}\in[x_{j,1}]$.
\end{lemma}

Define $\mathcal{O}_G$ as the digraph with the vertex set $G$,
and there is an arc from $x$ to $y$ if $y\prec x$.
Then $\mathcal{O}_G$ is an orientation
of $\mathcal P(G)$.

\begin{thm}\label{orientation}
Let $G$ be a group. Then $\mathcal{O}_G$ is a transitive orientation
of $\mathcal P_G$ and a subdigraph of $\overrightarrow{\mathcal{P}}_G$.
Moreover, if $\mathcal{O}'$ is a transitive orientation
of $\mathcal P_G$ and a subdigraph of $\overrightarrow{P}_G$,
then $\mathcal O'$ and $\mathcal O_G$ are isomorphic.
\end{thm}
\proof Suppose $\{(x,y),(y,z)\}\subseteq E(\mathcal{O}_G)$.
Then $z\prec y$ and $y\prec x$, and so$\langle z\rangle\subseteq\langle y\rangle\subseteq\langle x\rangle$.
If $\langle z\rangle\neq\langle x\rangle$, then $z\prec x$, which implies that $(x,z)\in E(\mathcal{O}_G)$.
If $\langle z\rangle=\langle x\rangle$, then $[z]=[y]=[x]$.
With reference to (\ref{1}), there exists an index $i$ such that
$z=x_{i,l}$, $y=x_{i,r}$, $x=x_{i,t}$ and $l<r<t$.  So $z\prec x$ and $(x,z)\in E(\mathcal{O}_G)$.
It follows that $\mathcal{O}_G$ is a transitive orientation of $\mathcal P_G$.
It is clear that $\mathcal{O}_G$ is  a subdigraph of $\overrightarrow{\mathcal{P}}_G$.

Assume that a subdigraph $\mathcal{O}'$ of $\overrightarrow{P}_G$
is another transitive orientation
of $\mathcal P_G$. With reference to (\ref{1}),
if $\langle x_{i,1}\rangle\subsetneqq\langle x_{j,1}\rangle$,
then $\{x_{i,1},x_{j,1}\}\in E(\mathcal P_G)$,
$(x_{i,1},x_{j,1})\in E(\overrightarrow{\mathcal P}_G)$ and
$(x_{j,1},x_{i,1})\not\in E(\overrightarrow{\mathcal P}_G)$,
which implies that $(x_{i,1},x_{j,1})\in E(\mathcal O')$.
Therefore, for $[x_{i,1}]\neq[x_{j,1}]$,
 we have $(x_{i,1},x_{j,1})\in E(\mathcal O_G)$ if and only
 if $(x_{i,1},x_{j,1})\in E(\mathcal O')$.
Since the induced subgraph on $[x_{i,1}]$
of $\mathcal P_G$ is a complete graph, and all transitive orientations
of a fixed complete graph are isomorphic,
we conclude that $\mathcal O'$ and $\mathcal O_G$ are isomorphic.
$\qed$

The following theorem is an immediate result from Theorem~\ref{orientation}.

\begin{thm}\label{comparability}
The power graph of a group is a
comparability graph.
\end{thm}

For two graphs $\Gamma$ and $\Gamma'$, a {\em homomorphism} from
$\Gamma$ to $\Gamma'$ is a
map $f:V(\Gamma)\longrightarrow V(\Gamma')$ such that
$\{f(u),f(v)\}\in E(\Gamma')$ whenever
$\{u,v\}\in E(\Gamma)$.
%A homomorphism from $\Gamma$ to the complete graph of order $k$ is {\em a $k$-coloring} of $\Gamma$.
The {\em chromatic number} of $\Gamma$,
denoted by $\chi(\Gamma)$, is the least value of $k$ such that
there exists a homomorphism from $\Gamma$ to the complete graph of order $k$.
%$\Gamma$ admits a $k$-coloring.
%A {\em clique} in $\Gamma$ is a subgraph of $\Gamma$ that is complete.
The {\em clique number} of $\Gamma$,
denoted by $\omega(\Gamma)$, is the maximum order of a clique in $\Gamma$.

A graph $\Gamma$ is {\em perfect} if $\chi(\Lambda)=\omega(\Lambda)$
for each induced subgraphs $\Lambda$ of $\Gamma$.
It was noted in \cite[Chapter V, Theorem~17]{bol} that comparability graphs
are perfect. Hence, by Theorem~\ref{comparability}, we have

\begin{cor}\label{pefect}{\rm\cite[Theorem 1]{doo}}
The power graph of a group is perfect.
\end{cor}

An {\em endomorphism} of a graph $\Gamma$ is a
homomorphism   from $\Gamma$ to itself.
A {\em core} \cite{god} of $\Gamma$ is a subgraph $\Lambda$
satisfies that every endomorphism of $\Lambda$
is an automorphism and there exists a
homomorphism from $\Gamma$ to $\Lambda$.
Every graph  has a core, which is an induced subgraph and is unique up to
isomorphism \cite[Lemma 6.2.2]{god}.
A graph is called a {\em core}
if its core is itself.
Godsil and Royle \cite{god2} showed that
the core of a graph $\Gamma$ is complete
if and only if $\chi(\Gamma)=\omega(\Gamma)$.

\begin{ob}\label{corcomplete}
The core of any induced subgraph of a perfect graph
is complete. In particular, the core of any induced
subgraph of a comparability graph is complete.
\end{ob}

\begin{prop}\label{complete}{\rm\cite[Theorem 2.12]{cha}}
Let $G$ be a group. Then $\mathcal{P}_G$ is complete
if and only if $G$ is a cyclic group of order $p^m$ for
some prime $p$ and nonnegative integer $m$.
\end{prop}

Combining Theorem~\ref{comparability}, Observation~\ref{corcomplete}
and Proposition~\ref{complete}, we get the following corollary.

\begin{cor}
Given a group $G$, the core of any induced subgraph of $\mathcal P_G$ is complete.
So $\mathcal P_G$ is a core if and only if $G$ is a cyclic group
of order $p^m$ for some prime $p$ and nonnegative integer $m$.
\end{cor}

\subsection{Posets}
A {\em partially ordered set} or {\em poset} $P$
is  an ordered pair $(V(P),\leq_P)$, where $V(P)$ is a finite set,
called the {\em vertex set} of $P$,
and $\leq_P$ is a reflexive, antisymmetric and transitive binary relation
on $V(P)$.
As usual, write $x<_P y$ if $x\leq_P y$ and $x\neq y$.
For any subset $S\subseteq V(P)$,
the {\em subposet} of $P$ induced by $S$, denoted by $P(S)$,
is a poset $(S,\leq_{P(S)}$), where $x\leq_{P(S)}y$ if and only if
$x\leq_P y$.
Two elements $x$ and $y$ of $V(P)$ are {\em comparable} if $x\leq_P y$
or $y\leq_P x$, otherwise $x$ and $y$ are {\em incomparable}.
The {\em comparability  graph}
of $P$, denoted by $\mathcal{G}_P$, is the graph with the
vertex set $V(P)$, where two distinct elements   are adjacent
if they are comparable.

From Theorem~\ref{orientation}, we get the following
example.
\begin{eg}\label{L(G)}
Let $G$ be a group. With reference to Definition~\ref{defin},
the ordered pair $(G,\preceq)$ is a poset.
In the remaining of this paper, we use $L_G$ to denote this poset.
The comparability graph of $L_G$ is the power graph of a group $G$, i.e., $\mathcal G_{L_G}=\mathcal P_G$.
\end{eg}

A {\em chain} (reps. An {\em antichain}) in a poset $P$ is
a subset of $V(P)$ such that all elements in this subset
are pairwise comparable (resp. incomparable).
A subset $S$ of $V(P)$ is {\em homogeneous}
if, for any
$y\in V(P)\setminus S$, one of the following holds:
\begin{itemize}
\item For all $x\in S$, $x\leq_P y$.
\item For all $x\in S$, $y\leq_P x$.
\item For all $x\in S$, $x$ and $y$ are incomparable.
\end{itemize}
A {\em homogeneous chain} (resp. {\em antichain}) in $P$
is a chain (resp. an antichain) that is homogeneous.
A partition $\mathcal S$ of $V(P)$ is
{\em a homogeneous partition} of $P$ if all elements
of $\mathcal S$ are homogeneous subsets.
Let $\mathcal S$ be a homogeneous partition of $P$.
The {\em quotient}  $P/\mathcal S=(\mathcal S,\leq_{P/\mathcal S})$, where two subsets $S_1,S_2\in\mathcal S$ satisfies $S_1\leq_{P/\mathcal S} S_2$ if $S_1=S_2$ or $x<_P y$ for each $x\in S_1$ and each $y\in S_2$. Then $P/\mathcal S$ is a poset.

The inverse operation of the
quotient is the {\em lexicographical sum} \cite{ille} defined as follows.
Let $P$ be a poset and let $\mathbb{Q}$ be a family of posets
indexed by $V(P)$, write $\mathbb{Q}=\{Q_x\mid x\in V(P)\}$.
The  lexicographical sum of $\mathbb Q$ over $P$,
denoted by $P[\mathbb Q]$,
is the poset with the vertex set $ V(P[\mathbb{Q}])=
\{(x,y)|x\in V(P)\textup{ and } y\in V(Q_x)\}$,
where $(x_1,y_1)\leq_{P[\mathbb{Q}]} (x_1,y_2)$ provided that
either $x_1=x_2$ and $y_1\leq_{Q_{x_1}}y_2$ or $x_1<_P x_2$.
 One can prove that this definition is well-defined. The following result is clear.

\begin{lemma}\label{inv}
 Suppose that $\mathcal S$ is a homogeneous partition of a poset $P$.
Write $R=P/\mathcal S$ and $\mathbb S=\{P(S)\mid S\in\mathcal S\}$.
Then $P$ is isomorphic to $R[\mathbb{S}]$.
\end{lemma}

Recall that the poset $L_G=(G,\preceq)$,
where $\preceq$ is defined in Definition~\ref{defin}.
The following lemma is an immediate result
from Lemma~\ref{lemma}.

\begin{lemma}\label{lemma quotient}
Let $G$ be a group.
With references to {\rm(\ref{1})},
any element $[x]$ in $\mathcal C'(G)$
is a homogeneous chain in $L_G$.
Consequently, the set $\mathcal C'(G)$ is a homogeneous partition
of $L_G$, and so $L_G/\mathcal C'(G)$ is a quotient of $L_G$.
\end{lemma}

The following result gives some equivalent conditions for comparing two distinct elements in the quotient $L_G/\mathcal C'(G)$.

\begin{lemma}\label{lemma1}
  Given a group $G$, let $[x]$ and $[y]$ be two distinct elements in $\mathcal C'(G)$. Then the following conditions are equivalent.

  {\rm(i)} $[x]<_{L_G/\mathcal C'(G)}[y]$.

  {\rm(ii)} $\langle x\rangle\subsetneqq\langle y\rangle$.

  {\rm(iii)} $(y,x)\in E(\overrightarrow{\mathcal P}_G)$ and $(x,y)\notin E(\overrightarrow{\mathcal P}_G)$
\end{lemma}
\proof If (i) holds, then $x\prec y$, which implies that (ii) holds by $[x]\neq[y]$. It is clear that (ii) and (iii) are equivalent. Suppose (ii) holds. Then, for each $x'\in[x]$ and each $y'\in[y]$, we have $\langle x'\rangle\subsetneqq\langle y'\rangle$, and so $x'\prec y'$. It follows that (i) holds.
$\qed$

\begin{prop}{\rm\cite[Theorem 2]{came}}\label{cam1}
If $G_1$ and $G_2$ are groups whose power graphs
are isomorphic, then their power digraphs are also isomorphic.
\end{prop}

%We use the following theorem to end this subsection.

\begin{thm}\label{lemma2}
Suppose $G_1$ and $G_2$ are groups. Then the followings are
equivalent.

{\rm(i)} The power graphs $\mathcal P_{G_1}$ and $\mathcal P_{G_2}$
are isomorphic.

{\rm(ii)} The  power digraphs $\overrightarrow{\mathcal P}_{G_1}$ and
$\overrightarrow{\mathcal P}_{G_2}$ are isomorphic.

{\rm(iii)} The transitive orientations $\mathcal O_{G_1}$
and $\mathcal O_{G_2}$ are isomorphic.

{\rm(iv)} The posets $L_{G_1}$ and $L_{G_2}$ are isomorphic.

{\rm(v)} There is an isomorphism $\tau$ from the quotient
$L_{G_1}/\mathcal C'(G_1)$ to the quotient $L_{G_2}/\mathcal C'(G_2)$ such that
$|\tau(S)|=|S|$ for each $S\in\mathcal C'(G_1)$.
\end{thm}
\proof Proposition~\ref{cam1} says that (i) implies
(ii). Lemma~\ref{orientation}
concludes that (ii) implies (iii).
From the definitions, we can see that (iii) and (iv) are equivalent.
By Example~\ref{L(G)},
it is clear that (iv) implies (i).
It follows from Lemmas~\ref{inv} and~\ref{lemma quotient}
that (v) implies (iv).

Suppose (ii) holds. Let $\sigma$ be an isomorphism
from  $\overrightarrow{\mathcal P}_{G_1}$ to
 $\overrightarrow{\mathcal P}_{G_2}$.
For $S\in\mathcal C'(G_1)$, define $\tau(S)=\{\sigma(x)\mid x\in S\}$.
Pick $x\in S$. Then
\begin{eqnarray*}
S=[x]&=&\{x\}\cup\{y\mid \{(x,y),(y,x)\}\subseteq E(\overrightarrow{\mathcal P}_{G_1})\}\\
&=&\{x\}\cup\{y\mid\{(\sigma(x),\sigma(y)),(\sigma(y),\sigma(x))\}\subseteq
E(\overrightarrow{\mathcal P}_{G_2})\},
\end{eqnarray*}
which implies that
$$
\tau([x])=\{\sigma(x)\}\cup \{\sigma(y)\mid\{(\sigma(x),\sigma(y)),(\sigma(y),\sigma(x))\}\subseteq
E(\overrightarrow{\mathcal P}_{G_2})\}=[\sigma(x)]\in\mathcal C'(G_2).
$$
Consequently, we obtain that $\tau$ is a bijection from
$\mathcal C'(G_1)$ to  $\mathcal C'(G_2)$ such that
$|\tau(S)|=|S|$. By Lemma~\ref{lemma1}, we have $S<_{L_{G_1}/\mathcal C'(G_1)}S'$  if and only if $\tau(S)<_{L_{G_2}/\mathcal C'(G_2)}\tau(S')$, and so  (v) holds.
$\qed$

\subsection{Characterization}
In order to give the structure of power graphs, we need the definition
of the generalized lexicographic product,
which was first defined by Sabidussi \cite{Sa}.
Given a graph $\mathcal{H}$ and a family of graphs
$\mathbb{F}=\{\mathcal{F}_v\mid v\in
V(\mathcal{H})\}$, indexed by $V(\mathcal{H})$, their {\em generalized lexicographic
product}, denoted by $\mathcal{H}[\mathbb{F}]$, is defined as the graph with
the vertex set $ V(\mathcal{H}[\mathbb{F}])=
\{(v,w)|v\in V(\mathcal{H})\textup{ and }
w\in V(\mathcal{F}_v)\} $ and the edge set $ E(\mathcal{H}[\mathbb{F}])
=\{\{(v_1,w_1),(v_2,w_2)\}|\{v_1, v_2\}\in E(\mathcal{H}), \textup{or }v_1=
v_2\textup{ and }\{w_1,w_2\}\in E(\mathcal{F}_{v_1})\}. $

Recall that the comparability graph of a poset $P$, denoted by $\mathcal G_P$, is the graph with the vertex set $V(P)$, where two distinct elements are adjacent if they are comparable.
\begin{lemma}\label{lex}
Given a poset $P$,  let $\mathbb{Q}$ be a family of posets indexed by $V(P)$.
Suppose $\mathbb{G}_{\mathbb{Q}}$ consists of all comparability graphs of posets in $\mathbb Q$.
Then $\mathcal G_{P[\mathbb Q]}=\mathcal G_P[\mathbb{G}_{\mathbb{Q}}]$.
\end{lemma}
\proof Write $\mathbb Q=\{Q_x\mid x\in V(P)\}$. It is clear that
$$
V(\mathcal G_{P[\mathbb Q]})=\{(x,y)\mid x\in V(P),y\in V(Q_x)\}=V(\mathcal G_P[\mathbb{G}_{\mathbb{Q}}]).
$$
Hence, it suffices to prove
$E(\mathcal G_{P[\mathbb Q]})=E(\mathcal G_P[\mathbb{G}_{\mathbb{Q}}])$.

Suppose $\{(x_1,y_1),(x_2,y_2)\}\in E(\mathcal G_{P[\mathbb Q]})$.
Then $(x_1,y_1)$ and $(x_2,y_2)$ are comparable in $P[\mathbb Q]$ and
$(x_1,y_1)\neq (x_2,y_2)$.
Without loss of generality, assume that
$(x_1,y_1)<_{P[\mathbb Q]} (x_2,y_2)$.
Hence, either $x_1=x_2$ and $y_1<_{Q_{x_1}}y_2$ or $x_1<_{P} x_2$.
If $x_1=x_2$ and $y_1<_{Q_{x_1}}y_2$,
then $\{y_1,y_2\}\in E(Q_{x_1})$, which implies that
$\{(x_1,y_1),(x_2,y_2)\}\in E(\mathcal G_P[\mathbb{G}_{\mathbb{Q}}])$.
If $x_1<_{P} x_2$, then $\{x_1,x_2\}\in E(\mathcal G_P)$. It follows that
$\{(x_1,y_1),(x_2,y_2)\}\in E(\mathcal G_P[\mathbb{G}_{\mathbb{Q}}])$.
Therefore, we have $E(\mathcal G_{P[\mathbb Q]})\subseteq E(\mathcal G_P[\mathbb{G}_{\mathbb{Q}}])$.

Suppose $\{(x_1,y_1),(x_2,y_2)\}\in E(\mathcal G_P[\mathbb{G}_{\mathbb{Q}}])$.
Hence, either $x_1=x_2$ and $\{y_1,y_2\}\in E(\mathcal G_{Q_{x_1}})$ or
$\{x_1,x_2\}\in E(\mathcal G_P)$.
If $x_1=x_2$ and $\{y_1,y_2\}\in E(\mathcal G_{Q_{x_1}})$, without loss of generality, assume that $y_1<_{Q_{x_1}}y_2$,
then $(x_1,y_1)<_{P[\mathbb{Q}]}(x_2,y_2)$, which implies that
$\{(x_1,y_1),(x_2,y_2)\}\in E(\mathcal G_{P[\mathbb Q]})$.
If $\{x_1,x_2\}\in E(\mathcal G_P)$,  assume that
$x_1 <_P x_2$, then $(x_1,y_1) <_{P[\mathbb{Q}]}(x_2,y_2)$, and so
$\{(x_1,y_1),(x_2,y_2)\}\in E(\mathcal G_{P[\mathbb Q]})$.  Therefore, we have $E(\mathcal G_P[\mathbb{G}_{\mathbb{Q}}])\subseteq E(\mathcal G_{P[\mathbb Q]})$.  We accomplish the proof.
$\qed$

Given a group $G$, let $\mathcal C(G)$ denote the set
of all cyclic subgroups of $G$.
Note that $(\mathcal C(G),\subseteq)$ is a poset. The following result is clear from Lemma~\ref{lemma1}.

\begin{lemma}\label{lemma3}
Let $G$ be a group. Then $L_G/\mathcal C'(G)$ is isomorphic to
 $(\mathcal C(G),\subseteq)$.
\end{lemma}

For a group $G$, define $\mathcal I_G$
as the graph with the vertex set $\mathcal C(G)$, and
two cyclic subgroups are adjacent if  one
is contained in the other. Then $\mathcal I_G$ is the comparability
graph of the poset $(\mathcal C(G),\subseteq)$.
For $C\in\mathcal{C}(G)$, let $\mathcal{K}_{C}$ be the complete graph of order
$\varphi(|C|)$, where $\varphi$ is the Euler's totient function.
Write $\mathbb{K}_G=\{\mathcal{K}_C\mid C\in\mathcal{C}(G)\}$.

\begin{thm}\label{str1}
Given a group $G$, the power graph $\mathcal{P}_G$ is isomorphic to
the generalized lexicographic product $\mathcal{I}_G[\mathbb{K}_G]$.
\end{thm}
\proof With reference to (\ref{1}) and by
 Definition~\ref{defin}, for any $i$, the subposet $L_G([x_{i,1}])$
is a totally ordered set, i.e., every pair of distinct elements  in $[x_{i,1}]$ are comparable. Therefore, the comparability
graph $\mathcal G_{L_G([x_{i,1}])}$ is the complete graph of
order $|[x_{i,1}]|$.
Since $|[x_{i,1}]|=\varphi(|\langle x_{i,1}\rangle|)$, we have
$\mathcal G_{L_G([x_{i,1}])}\simeq\mathcal{K}_{\langle x_{i,1}\rangle}$.
By Lemma~\ref{lemma3}, we obtain
$\mathcal G_{L_G/\mathcal C'(G)}\simeq\mathcal I_G$.
It follows that
\begin{equation}\label{2}
\mathcal G_{L_G/\mathcal C'(G)}[\{\mathcal G_{L_G([x_{i,1}])}\mid
[x_{i,1}]\in\mathcal C'(G)\}]\simeq \mathcal I_G[\mathbb{K}_G].
\end{equation}
By Lemma~\ref{inv} we get
$L_G\simeq (L_G/\mathcal C'(G))[\{L_G([x_{i,1}])\mid [x_{i,1}]\in\mathcal C'(G)\}]$.
Combining Example~\ref{L(G)}, Lemma~\ref{lex} and (\ref{2}), one has
$\mathcal P_G=\mathcal G_{L_G}\simeq\mathcal I_G[\mathbb{K}_G]$,
as desired.
$\qed$

Since the order of $\mathcal I_G[\mathbb{K}_G]$
is $\sum_{C\in\mathcal{C}(G)}\varphi(|C|)$, by Theorem~\ref{str1},
we get a Euler's classical formula on the finite groups.
If $G$ is cyclic, then this formula is the  Euler's classical formula.
\begin{cor}
Let $G$ be a group. Then
$$
\sum_{C\in\mathcal{C}(G)}\varphi(|C|)=|G|.
$$
\end{cor}

Finally, we give a necessary and sufficient condition for two isomorphic power graphs.

\begin{thm}
Let $G_1$ and $G_2$ be two groups. Then the followings
are equivalent.

{\rm(i)} The power graphs $\mathcal P_{G_1}$ and
$\mathcal P_{G_2}$ are isomorphic.

{\rm(ii)} There is an isomorphism $\sigma$ from
the poset $(\mathcal C(G_1),\subseteq)$ to the poset
$(\mathcal C(G_2),\subseteq)$ such that $|\sigma(C)|=|C|$
for each $C\in\mathcal C(G_1)$.
\end{thm}
\proof It follows from Theorem~\ref{lemma2} and Lemma~\ref{lemma3} that
(ii) implies (i).
Suppose (i) holds. By Theorem~\ref{lemma2} and Lemma~\ref{lemma3},
there exists an isomorphism $\sigma$ from $(\mathcal C(G_1),\subseteq)$ to
$(\mathcal C(G_2),\subseteq)$ such that $\varphi(|\sigma(C)|)=\varphi(|C|)$
for each $C\in\mathcal C(G_1)$, where $\varphi$ is the Euler's totient function.
In order to prove (ii), we only need to show that
$|\sigma(C)|=|C|$ for each $C\in\mathcal C(G_1)$.

Suppose for the contradiction that there exists
a cyclic subgroup $C_0$ of $G_1$ such that $|\sigma(C_0)|\neq|C_0|$.
Since $\varphi(|\sigma(C_0)|)=\varphi(|C_0|)$, there exists
a prime $p$ such that one of $|\sigma(C_0)|$
and $|C_0|$ is divided by $p$ and the other is not.

If $p$ divides $|C_0|$, there exists a cyclic subgraph
$C_1$ of $C_0$ with order $p$. Then
\begin{equation}\label{3}
\varphi(|\sigma(C_1)|)=\varphi(|C_1|)=p-1.
\end{equation}
For $i\in\{1,2\}$, let $e_i$ denote the identity of $G_i$. Then $\sigma(\langle e_1\rangle)=\langle e_2\rangle$.
Since there is no cyclic subgroup $C_2$ of $C_1$ such that
$\{e_1\}\subsetneqq C_2\subsetneqq C_1$, there is no cyclic subgroup
$C_2'$ of $\sigma(C_1)$ such that
$\{e_2\}\subsetneqq C_2'\subsetneqq\sigma(C_1)$,
which implies that $|\sigma(C_1)|$ is a prime,
and so $|\sigma(C_1)|=p$ by (\ref{3}).
Note that $\sigma(C_1)\subseteq\sigma(C_0)$.
Then $p$ divides $|\sigma(C_0)|$, a contradiction.
If $p$ divides $|\sigma(C_0)|$, we consider the inverse isomorphism
$\sigma^{-1}$, and similarly get a contradiction.
$\qed$

By the above theorem, we get the following proposition.

\begin{prop}{\rm\cite[Corollay 3]{came}}\label{cam2}
Two groups whose power graphs are isomorphic
have the same numbers of elements of each order.
\end{prop}

\section{Metric dimension}
In this section we establish a closed formula for the metric dimension
of the power graph of a group $G$. In Subsection 3.1, we give an equivalence relation on $G$ and denote by $\mathcal U(G)$ the set of all equivalence classes. If $G$ is cyclic, then $\mathcal U(G)$ is determined; otherwise we characterise all equivalence classes in $\mathcal U(G)$  by using  homogeneous sets  in $L_G$. In Subsection 3.2, we introduce a concept named resolving involution and denote by $W(G)$ the set of all resolving involutions of $G$. If $G$ is cyclic, then $W(G)$ is determined; otherwise, by using homogeneous sets in a subposet of $L_G$, we provide a necessary and sufficient condition for an involution to be a resolving involution of $G$. In Subsection 3.3, we establish a closed formula for $\dim(\mathcal P_G)$ in terms of  $|G|$, $|\mathcal U(G)|$ and $|W(G)|$. In particular, we compute the metric dimension of the power graph of a cyclic group.

\subsection{Equivalence classes}
Given an element $x$ in a group $G$, the {\em open neighborhood} of $x$ in the power graph $\mathcal P_G$, denoted by  $N(x)$, is the set $\{y\in G\mid d_{\mathcal P_G}(x,y)=1\}$;
the {\em closed neighborhood} of $x$ in $\mathcal P_G$, denoted by $N[x]$, is the union of $N(x)$ and $\{x\}$.

For two elements $x$ and $y$ in a group $G$,
define $x\equiv y$ if $N(x)=N(y)$ or $N[x]=N[y]$. Hernando et al. \cite{He}
proved that $\equiv$ is an equivalence relation.
Let $\overline{x}$ denote the equivalence class that contains $x$.
%Then $\overline x$ is a clique or an independent set in $\mathcal P_G$.
Write
$$
\mathcal U(G)=\{\overline{x}\mid x\in G\}.
$$

\begin{ob}\label{equi}
Let $x$ be an element of a group $G$.

(i)  $[x]\subseteq\overline x$.

(ii) $\overline x=\{y\mid y\in G, N(y)=N(x)\}$ or
$\{y\mid y\in G, N[y]=N[x]\}$. In particular, the equivalence class $\overline x$ is
an independent set or a clique in $\mathcal P_G$.
\end{ob}

A {\em maximal involution} of a group $G$ is an involution $x$ such that
$\langle x\rangle$ is a maximal cyclic subgroup of $G$. For $y\in G$, let $o(y)$ denote the order of $y$ in the rest of this paper.

\begin{lemma}\label{open}
Let $G$ be a group of order at least two.
Suppose that $x$ and $y$ are two distinct elements in $G$.
Then $N(x)=N(y)$ if and only if both $x$ and $y$ are
maximal involutions of $G$.
\end{lemma}
\proof If both $x$ and $y$ are
maximal involutions, then $G$ is noncyclic, and so
$N(x)=\{e\}=N(y)$.
Now suppose $N(x)=N(y)$.

%By Observation~\ref{equi}, we have $[x]\subseteq\overline x=\overline y$
%and $\overline y$ is an independent set in $\mathcal P_G$,
%which implies that $[x]$ is an independent set, and so $[x]=\{x\}$.
%Hence, we get $o(x)=2$. Similarly, we have $o(y)=2$.

If $o(x)\geq 3$, then $x^{-1}\neq x$.
Note that $N[x]=N[x^{-1}]$.
Since $x^{-1}\in N(x)$, we have $x^{-1}\in N(y)$,
and so $y\in N(x^{-1})$, which implies that
$y\in N(x)$, a contradiction.
So $o(x)=2$. Similarly, we have $o(y)=2$.

If $\langle x\rangle$ is not a maximal cyclic subgroup,
there exists an element
$z$ of even order in $G\setminus\{x\}$ such that
$\langle x\rangle\subseteq\langle z\rangle$, which implies that
$z\in N(x)$, and so $z\in N(y)$.
Consequently, one gets $\langle y\rangle\subseteq\langle z\rangle$.
Note that the involution in a cyclic group of
even order is unique. Hence $x=y$, a contradiction.
Therefore $\langle x\rangle$ is a maximal cyclic subgroup.
We obtain that $\langle y\rangle$ is a maximal cyclic subgroup similarly.
$\qed$

\begin{lemma}\label{subequiv}
Given a group $G$, let $U$ be a homogeneous antichain or a
homogeneous chain in $L_G$. Then $U\subseteq\overline u$ for any $u\in U$.
\end{lemma}
\proof Pick $x,y\in U$ and $z\in G\setminus U$. Since $U$ is homogeneous,
we have $z\in N(x)$ is equivalent to $z\in N(y)$. Hence,
if $U$ is an antichain, then $N(x)=N(y)$; if $U$ is a chain, then
$N[x]=N[y]$. Consequently, the desired result follows.
$\qed$

Let $G$ be a group. The following result, the proof of which is immediate from  Lemmas~\ref{open} and \ref{subequiv}, characterise equivalence classes in $\mathcal U(G)$ that is an independent set with at least two vertices in $\mathcal P_G$.

\begin{prop}\label{homanti}
Suppose that $U$ is a subset of a group $G$ and  $|U|\geq 2$.
Then the followings are equivalent.

{\rm(i)} The set $U$ is an equivalence class in $\mathcal U(G)$
that is an independent set in $\mathcal P_G$.

{\rm(ii)} The set $U$ consists of all maximal involutions of $G$.

{\rm(iii)} The set $U$ is a maximal homogeneous antichain in $L_G$.
\end{prop}

Given a group $G$, we always use $e$ to denote the identity in the remaining of this paper. Now we consider the equivalence class in $\mathcal U(G)$ that is
a clique in $\mathcal P_G$. Note that $\overline e$ is always a clique in $\mathcal P_G$.

\begin{prop}\label{identity}{\rm\cite[Proposition 4]{cam}}
{\rm(i)} Suppose $G$ is a cyclic group generated by $x$.
If $|G|$ is a prime power, then $\overline e=G$.
If $|G|$ is not a prime power,
then $\overline e=[e]\cup[x]$.

{\rm(ii)} If $G$ is a generalized quaternion  $2$-group, then $\overline e=\{e,x\}$,
where $x$ is the unique
involution in $G$.

{\rm(iii)} If a group $G$ is neither a cyclic group nor
a generalized quaternion  $2$-group, then $\overline e=\{e\}$.
\end{prop}

\begin{prop}\label{closed}{\rm\cite[Proposition 5]{cam}}
Let $x$ be an element of a group $G$. Suppose that $\overline x\neq \overline e$ and $\overline x$  is
a clique in $\mathcal P_G$. Then
one of the following holds.

{\rm(i)} $\overline x=[x]$.

{\rm(ii)} There exist elements $x_0,x_1,\ldots,x_r$ in $G$ with
$\langle x_0\rangle\subsetneqq\langle x_1\rangle\subsetneqq\cdots\subsetneqq\langle x_r\rangle$
and $o(x_i)=p^{s+i}$ for each $i\in\{0,1,\ldots,r\}$, where $p$ is a prime and $s$ is
a positive integer, such that
$$
\overline x=[x_0]\cup[x_1]\cup\cdots\cup[x_r].
$$
\end{prop}

For a cyclic group $G$, the set $\mathcal U(G)$ is determined in the following result.

\begin{prop}{\label{cycequi}}
Suppose a cyclic group  $G=\langle x\rangle$.

{\rm(i)} If $o(x)$ is a prime
power, then $\mathcal U(G)=\{G\}$.

{\rm(ii)} If $o(x)$ is not a prime power, then
$$
\mathcal U(G)=(\mathcal C'(G)\setminus\{[e],[x]\})\cup\{[e]\cup[x]\}.
$$
\end{prop}
\proof (i) It is immediate from Proposition~\ref{identity}.

(ii) By Proposition~\ref{identity}, we only need to show that $\overline y=[y]$
for any $y\in G\setminus([e]\cup[x])$.
If $\overline y\neq[y]$, by Proposition~\ref{closed}, we have $o(y)=p^s$ and
there exists an element
$z\in G\setminus[y]$ with $o(z)=p^t$ such that $[y]\cup[z]\subseteq\overline y$,
where $p$ is a prime and $s,t$ are two distinct positive integers.
Choose a prime divisor $q$ of $|G|$ with $q\neq p$. Hence,
there exists an element $u\in G$ with $o(u)=p^{\max\{s,t\}}q$, which implies that $u\in N[y]\setminus N[z]$ or $u\in N[z]\setminus N[y]$,
and so $N[y]\neq N[z]$, a contradiction.
$\qed$

\begin{lemma}\label{lemma4}
Let $x$ be an element of a noncyclic group $G$. If $\overline x$ is a clique in $\mathcal P_G$, then $\overline x$ is a homogeneous chain in $L_G$.
\end{lemma}
\proof If $\overline x=\overline e$, then $G$ is a generalized quaternion
$2$-group and $\overline x=\{e,x_0\}$ by Propositions~\ref{identity},
where $x_0$ is the unique involution.
For any $y\in G\setminus\overline x$, we have $x_0=y^{\frac{o(y)}{2}}$
and $e=y^{o(y)}$, which implies that $x_0\prec y$
and $e\prec y$. So $\overline x$ is a homogeneous chain in $L_G$.
Now suppose $\overline x\neq\overline e$. Then (i) or (ii) in Proposition~\ref{closed} holds.
If (i) holds, then $\overline x$ is a homogeneous chain in $L_G$ by Lemma~\ref{lemma quotient}.
Suppose (ii) holds. Then   there exist elements  $x_1$ and $x_2$ in $\overline x$ with $o(x_1)=p^s$ and $o(x_2)=p^t$, where $p$ is a prime and $s<t$, such that
$\overline x=\{y\mid \langle x_1\rangle\subseteq\langle y\rangle\subseteq\langle x_2\rangle\}$.
If $\overline x$ is not a homogeneous chain in $L_G$, there exist elements $z\in G\setminus \overline x$ and $y_1,y_2\in\overline x$ such that
$y_1\prec z\prec y_2$,
then
$\langle x_1\rangle\subseteq\langle y_1\rangle\subsetneqq \langle z\rangle\subsetneqq \langle y_2\rangle\subseteq\langle x_2\rangle$,
and so $z\in\overline x$, a contradiction.
$\qed$

Let $G$ be a noncyclic group. In  the following two propositions, by using  homogeneous sets  in $L_G$, we characterise equivalence classes in $\mathcal U(G)$ that is a clique in $\mathcal P_G$. The proof of Proposition~\ref{homchai}  is clear from Lemmas~\ref{subequiv} and \ref{lemma4}, and the proof of Proposition~\ref{one} is immediate from Propositions~\ref{homanti} and \ref{homchai}.

\begin{prop}\label{homchai}
Suppose that $U$ is a subset of a noncyclic group $G$ and $|U|\geq 2$. Then $U$ is an
equivalence class in $\mathcal U(G)$ that is
a clique in $\mathcal P_G$ if and only if $U$
is a maximal homogeneous chain in $L_G$.
\end{prop}

\begin{prop}\label{one}
  Let $x$ be an element of a noncyclic group $G$. Then $\{x\}\in\mathcal U(G)$ if and only if $\{x\}$ is a maximal homogeneous chain and a maximal homogeneous antichain in $L_G$.
\end{prop}

\subsection{Resolving involutions}
We begin this subsection by a notation. For elements $x$ and $y$ in a group $G$, write
$$
R\{x,y\}=\{z\mid z\in G, d_{\mathcal P_G}(x,z)\neq d_{\mathcal P_G}(y,z)\}.
$$

\begin{ob}\label{rxy}
Let $G$ be a group. Pick two distinct elements $x$ and $y$.

(i) Any resolving set of $\mathcal P_G$
intersects $R\{x,y\}$ nonempty.

(ii) The equation $\overline x=\overline y$ holds if and only if $R\{x,y\}=\{x,y\}$.

(iii) If there exists an element $z\in R\{x,y\}\setminus\{x,y\}$,  then $\overline z\subseteq R\{x,y\}$.
\end{ob}

A {\em resolving involution} of a group $G$ is
an involution $w$ satisfies that there exist two elements $x,y\in G\setminus\overline w$ with $R\{x,y\}=\{x,y,w\}$.
%In this subsection, we always assume that $G$ has a resolving involution.
Let $W(G)$ denote the set of all resolving involutions of $G$.
For each $w\in W(G)$, fix two elements $x_w$ and $y_w$
such that $R\{x_w,y_w\}=\{x_w,y_w,w\}$.

\begin{ob}\label{ri}
Suppose that $w$ is a resolving involution of a group $G$.

(i) Then $\overline w=\{w\}$.

(ii) Then $\overline{x_w},\overline{y_w}$ and $\overline w$ are pairwise
distinct.

(iii) For each pair
$(x,y)\in\overline{x_w}\times\overline{y_w}$, we have
$R\{x,y\}=\{x,y,w\}$.
\end{ob}

\begin{lemma}\label{e2}
Let $w$ be a resolving involution of a group $G$.
Then $\langle x_w\rangle\subseteq\langle y_w\rangle$
or $\langle y_w\rangle\subseteq\langle x_w\rangle$.
\end{lemma}
\proof Suppose for the contrary that
$\langle x_w\rangle\nsubseteq\langle y_w\rangle$
 and $\langle y_w\rangle\nsubseteq\langle x_w\rangle$.
Then $x_w$ and $y_w$ are not adjacent in $\mathcal P_G$.
Hence $[x_w]\cup[y_w]\subseteq R\{x_w,y_w\}=\{x_w,y_w,w\}$.
By Observations~\ref{equi} and \ref{ri}, we have
$[x_w]=\{x_w\}$ and $[y_w]=\{y_w\}$, which implies that $o(x_w)=o(y_w)=2$,
and so $w\not\in R\{x_w,y_w\}$, a contradiction.
$\qed$

\begin{lemma}\label{ril1}
Let $w$ be a resolving involution of a group $G$. Then there exists
a cyclic subgroup $C$ of $G$ such that
$w$ is the resolving involution of $C$.
\end{lemma}
\proof  It suffices to show that
there exists a cyclic subgroup $C$ of $G$ such that $\{x_w,y_w,w\}\subseteq C$.
By Lemma~\ref{e2}, without loss of generality, assume that $\langle x_w\rangle\subseteq\langle y_w\rangle$.
Hence, we only need to consider $w\not\in\langle y_w\rangle$.
Since $w\in R\{x_w,y_w\}$, we have $x_w\in\langle w\rangle$,
which implies that $x_w=e$.

{\em Claim 1. For any $z\in G\setminus\{w\}$,
we have $\langle z\rangle\subseteq\langle y_w\rangle$
or $\langle y_w\rangle\subsetneqq\langle z\rangle$.}
In fact, if $z\in G\setminus\{e,y_w,w\}$, then $z\not\in R\{e,y_w\}$,
which implies that $z$ is adjacent to $y_w$ in $\mathcal P_G$.
Hence, Claim 1 is valid.

Write
$$
\mathcal A=\{\langle z\rangle\mid z\in G\setminus\{w\},\langle y_w\rangle\subsetneqq\langle z\rangle\}.
$$
If $\mathcal A=\emptyset$, by Claim 1, we have
$\langle y_w\rangle=G\setminus\{w\}$, which implies that
$(|G|-1)$ is a divisor of $|G|$, a contradiction.
So $\mathcal A\neq\emptyset$.

{\em Claim 2. For any $\langle z\rangle\in\mathcal A$,
if $w\not\in\langle z\rangle$, then $o(z)$ is a prime
power.}
In fact, if $o(z)$ is not a prime power, since $o(y_w)$
divides $o(z)$ and $o(y_w)\neq o(z)$, there exists a divisor $m$
of $o(z)$ such that $m$ does not divide $o(y_w)$ and
$o(y_w)$ does not divide $m$. Pick $z_0\in\langle z\rangle$ with $o(z_0)=m$.
Then $z_0\neq w$, $\langle z_0\rangle\nsubseteq \langle y_w\rangle$ and
$\langle y_w\rangle\nsubseteq \langle z_0\rangle$,  contrary to Claim 1.
Hence, Claim 2 holds.

Suppose that $w\not\in\langle z\rangle$ for any $\langle z\rangle\in\mathcal A$.
 By Claim 2, it is clear that $o(y_w)$ is a prime power.
Write $o(y_w)=p^s$, where $p$ is a prime and $s$ is a positive integer.
By Claims 1 and 2, the following claim is valid.

{\em Claim 3. For any $z\in G\setminus\{w\}$,
we get $o(z)=p^i$ for some nonnegative integer $i$.}

{\em Claim 4. The subgroup of order $p$ that is contained in $G\setminus\{w\}$
is unique.} In fact, the subgroup of order $p$  in $\langle y_w\rangle$
is unique, which we denote by $P$. If there exists two subgroups $P$ and $Q$
of order $p$ such  that $P\cup Q\subseteq G\setminus\{w\}$,
then $Q\cap\langle y_w\rangle=\{e\}$, contrary to Claim 1.
Hence, Claim 4 holds.

Write $m_i=|\{x\mid x\in G\setminus\{w\},o(x)=p^i\}|$.
Let $t$ be the maximum number of $i$ such that $m_i\neq 0$.
By Claim 3, we have
\begin{equation}\label{e1}
\sum_{i=0}^t m_i=|G|-1.
\end{equation}
Since  $\varphi(p^i)$ divides $m_i$,
the prime $p$ divides $m_i$ for $i\in\{2,\ldots,t\}$,
which implies that $p$ divides $|G|-1-m_0-m_1$ by (\ref{e1}).
It is clear that $p$ divides $|G|$ and $m_0=1$.
So $p$ divides $m_1+2$.
By Claim 4, we have $m_1=p-1$, which implies that
$p$ divides $p+1$, a contradiction.

Therefore, there exists a cyclic subgroup $\langle z\rangle\in\mathcal A$
with $w\in\langle z\rangle$, which implies that $\{x_w, y_w,w\}\subseteq\langle z\rangle$.
$\qed$

For a cyclic group $G$, the set $W(G)$ is determined in the following result.

\begin{prop}\label{ril2}
Let $w$ be the involution of a cyclic group $G$.

{\rm(i)} If $w$ is a resolving involution, then $|G|=2p^m$ or
$2^mp$ for some positive integer $m$ and odd prime $p$.

{\rm(ii)} If
$|G|=2p^m$, then $w$ is a resolving involution and
$$
\{o(x_w),o(y_w)\}\in\{\{1,p\},\{2p^m,p\}\}.
$$

{\rm(iii)} If
$|G|=2^m p$ and $m\geq 2$, then $w$ is a resolving involution and
$$
\{o(x_w),o(y_w)\}=\{2p, p\}.
$$

\end{prop}
\proof (i)  Write $G'=G\setminus\{x_w,y_w,w\}$. Since $R\{x_w,y_w\}=\{x_w,y_w,w\}$, we get the following claim.

{\em Claim 1. For any $z\in G'$, we have $d_{P_G}(z,x_w)=d_{P_G}(z,y_w)$. }

Since $w\in R\{x_w,y_w\}$,
in $\mathcal P_G$ one of $x_w$ and $y_w$ is adjacent to $w$ and the other is not.
Without loss of generality, assume that $x_w$ and $w$ are adjacent.
Then $y_w$ and $w$ are not adjacent.
Hence , the following claim is valid.

{\em Claim 2. The number $o(y_w)$ is  odd  and $o(y_w)\geq 3$.}

Write $|G|=2^{s_0}p_1^{s_1}\cdots p_t^{s_t}$,
where $p_1,\ldots,p_t$ are odd primes and $s_0,s_1,\ldots,s_t$
are positive integers.
Now we divide our proof into two cases.

{\em Case 1.} $\overline{x_w}=\overline e$. Then $x_w$ is adjacent to any element of $G'$ in $\mathcal P_G$.
If $s_0\geq 2$,  there exists an element $z_0\in G'$ of order $4$, then $z_0$ and $y_w$ are adjacent by Claim 1, and so $4$ divides $o(y_w)$ or $o(y_w)$ divides $4$, contrary to Claim 2. Hence $s_0=1$.
If $t\geq 2$, then there exist elements $z_1$ and $z_2$ in $G'$ of order $2p_1$ and $2p_2$, respectively. By Claim 1, both $z_1$ and $z_2$ are adjacent to $y_w$. It follows from Claim 2 that $o(y_w)=p_1=p_2$, a contradiction.
Hence $t=1$ and $o(y_w)=p_1$.
So $|G|=2p_1^{s_1}$.

{\em Case 2.} $\overline{x_w}\neq\overline e$.
Then  $o(x_w)$ is even. Write
$o(x_w)=2^{i_0}p_{j_1}^{i_1}\cdots p_{j_l}^{i_l}$,
where $\{j_1,\ldots,j_l\}\subseteq\{1,\ldots,t\}$
and $1\leq i_k\leq s_k$ for $k\in\{0,1,\ldots,l\}$.
Similar to Case 1, we get $i_0=1$, $l=1$ and $o(y_w)=p_{j_1}$.
So $o(x_w)=2p_{j_1}^{i_1}$.
If $t\geq 2$,   there exists an elements $z_3\in G'$
such that $o(z_3)=p_{j_1}q$ for some prime
$q\in\{p_1,\ldots,p_t\}\setminus\{p_{j_1}\}$, then
$z_3$ is adjacent to $y_w$ and not adjacent to $x_w$,
contrary to Claim 1. Hence $t=1$.
Consequently, we get $|G|=2^{s_0}p_1^{s_1}$,
$o(x_w)=2 p_1^{i_1}$ and $o(y_w)=p_1$.
If $i_1<s_1$, then any element of order
$p_1^{i_1+1}$ in $G'$ is adjacent to $y_w$ and not adjacent to $x_w$,
contrary to Claim 1.
Therefore $o(x_w)=2p_1^{s_1}$.
If $s_1\geq 2$,
then any element of order
$p_1^2$ in $G'$ is adjacent to $y_w$ and not adjacent to $x_w$,
contrary to Claim 1. Hence $s_1=1$, and so$|G|=2^{s_0}p_1$.

(ii) Suppose $y$ is an element of $G$ with
$o(y)=p$. Then $R\{e,y\}=\{e,y,w\}$, which implies that
$w$ is a resolving involution of $G$. Combining Proposition~\ref{identity} and the proof of  (i), we have $\{o(x_w),o(y_w)\}\in\{\{1,p\},\{2p^m,p\}\}.$

(iii) Suppose that $x_1$ and $x_2$ are two elements of $G$ with
$o(x_1)=2p$ and $o(x_2)=p$.
Then $R\{x_1,x_2\}=\{x_1,x_2,w\}$, which implies that
$w$ is a resolving involution of $G$.
It follows from the proof of (i) that
$\{o(x_w),o(y_w)\}=\{2p, p\}.$
$\qed$

\begin{lemma}{\rm\cite[Theorem 5.4.10. (ii)]{gor}}\label{p}
Let $p$ be a prime. If $G$ is a $p$-group which has a unique minimal
subgroup of order $p$, then $G$ is either a cyclic group or a generalized quaternion group.
\end{lemma}

In the rest of this subsection, we consider  the resolving involutions of a noncyclic group.

\begin{prop}\label{e6}
Let $w$ be a resolving involution of a noncyclic group $G$. Suppose $o(x_w)\leq o(y_w)$. Then the follows hold.

{\rm(i)} $\langle x_w\rangle\cup\langle w\rangle\subseteq\langle y_w\rangle$.

{\rm(ii)} There exists an odd prime divisor $p$ of $|G|$ such that $(o(x_w),o(y_w))=(p,2p^m)$ for some positive integer $m$.
\end{prop}
\proof   In order to prove (i) and (ii), combining Lemma~\ref{ril1} and Proposition~\ref{ril2}, we only need to show that $(o(x_w),o(y_w))\neq(1,q)$ for any odd prime $q$. Suppose for the contrary that $x_w=e$ and $o(y_w)=q$ for some odd prime $q$. Since $R\{e,y_w\}=\{e,y_w,w\}$, each element in $G\setminus\{w,y_w\}$ is adjacent to $y_w$ in $\mathcal P_G$, which implies that $y_w\in\langle z\rangle$ for any $z\in G\setminus\{e,w\}$. Hence, the following claims are valid.

{\em Claim 1.  All prime divisors of $|G|$ are $2$ and $q$.}

{\em Claim 2. The group $G$ contains a unique involution, which is  $w$, and a unique subgroup of order $q$, which is $\langle y_w\rangle$.}

{\em Claim 3. There is no element of order $4$ in $G$. }

By Claims 2 and 3,   the subgroup $\langle w\rangle$ is a unique Sylow $2$-subgroup of $G$, and so $\langle w\rangle$ is normal in $G$. By Claim 1, we have $|G|=2q^n$ for some positive integer $n$. By Claim 2 and Lemma~\ref{p},  a Sylow $q$-subgroup $Q$ of $G$  is isomorphic to the cyclic group of order $q^n$.
Since the index of $Q$ in $G$ is $2$,  the Sylow $q$-subgroup $Q$ is normal in $G$. Consequently, the group $G$ is isomorphic to $\langle w\rangle\times Q$, which is isomorphic to the cyclic group of order $2q^n$,  a contradiction.
$\qed$

Given a noncyclic group $G$, by using the homogeneous set  in a subposet of $L_G$, we provide a necessary and sufficient condition for an involution to be a resolving involution of $G$.

\begin{prop}
Let $w$ be an involution of a noncyclic group $G$. Then $w$ is a resolving involution of $G$ if and only if there exists a cyclic subgroup $C$  of $G$ such that the following conditions hold.

 {\rm (i)}  $|C|=2p^m$ for some odd prime $p$ and positive integer $m$.

{\rm (ii)} $w\in C$.

{\rm (iii)} The set $C\setminus\langle w\rangle$ is homogeneous in the subposet $L_G(G\setminus\{w\})$.
\end{prop}
\proof Suppose $w$ is a resolving involution of $G$. Without loss of generality, assume that $o(x_w)\leq o(y_w)$. By Proposition~\ref{e6}, we have $\langle x_w\rangle\cup\langle w\rangle\subseteq\langle y_w\rangle$ and $(o(x_w),o(y_w))=(p,2p^m)$ for some odd prime $p$ and positive integer $m$. Let $C=\langle y_w\rangle$. Then (i) and (ii) hold. Now we prove (iii).

For each $x\in C\setminus\langle w\rangle$, since $\langle x_w\rangle\subseteq\langle x\rangle\subseteq\langle y_w\rangle$, we have
\begin{equation}\label{8}
  x_w\preceq x\preceq y_w.
\end{equation}
Pick any $z\in (G\setminus\{w\})\setminus(C\setminus\langle w\rangle)$. Then $z=e$ or $z\in G\setminus C$. If $z=e$, then $z\preceq x$ for each $x\in C\setminus\langle w\rangle$. In the following two cases, suppose $z\in G\setminus C$.

{\em Case 1.} $C\subseteq\langle z\rangle$. Then $y_w\preceq z$. For each $x\in C\setminus\langle w\rangle$, by (\ref{8}), we have  $x\preceq z$.

{\em Case 2.} $C\nsubseteq\langle z\rangle$. If there exists an element $x_1\in C\setminus\langle w\rangle$ such that $z\preceq x_1$, then $z\preceq y_w$ by (\ref{8}). So $z\in C$, a contradiction.  If there exists an element $x_2\in C\setminus\langle w\rangle$ such that $x_2\preceq z$, then $x_w\preceq z$ by (\ref{8}).  Hence $z$ is adjacent to $x_w$ in $\mathcal P_G$. Since $z\not\in\{x_w,y_w,w\}=R\{x_w,y_w\}$, elements $z$ and $y_w$ are adjacent, which implies that $z\in C$ or $C\subseteq\langle z\rangle$, a contradiction. Therefore, for each $x\in C\setminus\langle w\rangle$, elements $z$ and $x$ are incomparable in $L_G(G\setminus\{w\})$.

Hence (iii) holds.

Conversely, if there exists a cyclic subgroup $C$ of $G$ such that (i), (ii) and (iii) hold. Write $C=\langle y\rangle$. By (i) and (ii), we get $o(y)=2p^m$ and $w=y^{p^m}$, where $p$ is an odd prime and $m$ is a positive integer. Consider these two vertices $y$ and $y^{2p^{m-1}}$ in $P_G$.
By (iii), any vertex in $ (G\setminus\{w\})\setminus(\langle y\rangle\setminus\langle w\rangle)$ is adjacent to both or neither of  them.
Note that each vertex in $\langle y\rangle\setminus\{w,y,y^{2p^{m-1}}\}$ is adjacent to both of them. Hence, we have $R\{y,y^{2p^{m-1}}\}=\{w,y,y^{2p^{m-1}}\}$. It follows that $w$ is a resolving involution of $G$.
$\qed$

The following lemma is useful for the next subsection.

\begin{lemma}\label{ril3}
Let $u$ and $v$ be two distinct resolving involutions of a group $G$.
Then
$\{\overline{x_u},\overline{y_u},\overline{u}\}
\cap\{\overline{x_v},\overline{y_v},\overline{v}\}=\emptyset$.
\end{lemma}
\proof  Without loss of generality, assume that $o(x_u)\leq o(y_u)$ and $o(x_v)\leq o(y_v)$.
Since $G$ has at least two involutions, we know that $G$ is noncyclic.  By Proposition~\ref{e6},
we get
\begin{equation*}
\langle x_w\rangle\cup\langle w\rangle\subseteq\langle y_w\rangle\quad\textup{and}\quad(o(x_w),o(y_w))\in\{(p,2p^m)\mid p \textup{ is an odd prime}, m\geq 1\},
\end{equation*}
where $w\in\{u,v\}$. Then $\langle y_u\rangle\nsubseteq\langle y_v\rangle$ and $\langle y_v\rangle\nsubseteq\langle y_u\rangle$. So $\overline{y_u}\neq\overline{y_v}$.
If  $\overline{x_u}=\overline{x_v}$, then $y_v$ is adjacent to $x_u$ in $\mathcal P_G$, which implies that $y_v\in R\{x_u,y_u\}$, a contradiction. Hence $\overline{x_u}\neq\overline{x_v}$. It follows that
$\{\overline{x_u},\overline{y_u},\overline{u}\}
\cap\{\overline{x_v},\overline{y_v},\overline{v}\}=\emptyset$, as desired.
$\qed$

\subsection{Formula}
In this subsection, we shall establish a closed formula for the metric dimension of the power graph of a group. As an application, we compute $\dim(\mathcal P_{Z_n})$, where $Z_n$ is a cyclic group of order $n$. We begin by some lemmas.

\begin{lemma}\label{e7}
Let $G$ be a group. Suppose that $S$ is a resolving set of $\mathcal P_G$ and $\overline z\in\mathcal U(G)$. Then $|S\cap\overline z|\geq |\overline z|-1$.
\end{lemma}
\proof If $|S\cap\overline z|\leq|\overline z|-1$, there exist  two distinct elements $z_1,z_2\in\overline z$ such that $S\cap\{z_1,z_2\}=\emptyset$. Since $\overline{z_1}=\overline{z_2}=\overline z$, by Observation~\ref{rxy} we have $R\{z_1,z_2\}=\{z_1,z_2\}$ and $S\cap R\{z_1,z_2\}\neq\emptyset$, a contradiction.
$\qed$

\begin{lemma}\label{low}
Given a group $G$, we have
 $\dim(\mathcal P_G)\geq|G|-|\mathcal U(G)|+|W(G)|$.
\end{lemma}
\proof Suppose that $S$ is a resolving set of $\mathcal P_G$ with
size $\dim(\mathcal P_G)$. If $W(G)=\emptyset$, by Lemma~\ref{e7} we get
$$
\dim(\mathcal P_G)=|S|=\sum_{\overline z\in\mathcal U(G)}|S\cap\overline z|\geq \sum_{\overline z\in\mathcal U(G)}(|\overline z|-1)=|G|-|\mathcal U(G)|.
$$

Now suppose $W(G)\neq\emptyset$.
For each $w\in W(G)$, by Observations~\ref{ri} and Lemma~\ref{e7}, we get
\begin{equation*}
|S\cap(\overline{w}\cup\overline{x_w}\cup\overline{y_w})|\geq |\overline w|-1+|\overline{x_w}|-1+|\overline{y_w}|-1+1=
|\overline{x_w}|+|\overline{y_w}|-1,
\end{equation*}
which implies that
\begin{eqnarray}\label{6}
\sum_{w\in W(G)}(|S\cap\overline w|+|S\cap\overline{x_{w}}|+|S\cap\overline{y_{w}}|)
\geq\sum_{w\in W(G)}(|\overline{x_{w}}|+|\overline{y_{w}}|)-|W(G)|.
\end{eqnarray}
Write
$
\mathcal W(G)=\bigcup_{w\in W(G)}\{\overline{w},\overline{x_{w}},\overline{y_{w}}\}.
$
Combining Lemma~\ref{ril3} and (\ref{6}), we have
\begin{eqnarray}
\sum_{\overline z\in\mathcal W(G)}|S\cap\overline z|
%&=&\sum_{w\in W(G)}(|S\cap\overline{w}|+|S\cap\overline{x_{w}}|+|S\cap\overline{y_{w}}|)\nonumber\\
\geq\sum_{w\in W(G)}(|\overline{x_{w}}|+|\overline{y_{w}}|)-|W(G)|
=\sum_{\overline z\in\mathcal W(G)}|\overline z|-2|W(G)|.\label{7}
\end{eqnarray}
By (\ref{7}) and Lemma~\ref{e7}, we get
\begin{eqnarray*}
\dim(\mathcal P_G)=|S|
&=&\sum_{\overline z\in\mathcal W(G)}|S\cap\overline z|+
\sum_{\overline z\in\mathcal U(G)\setminus\mathcal W(G)}|S\cap\overline z|\\
&\geq&\sum_{\overline z\in\mathcal W(G)}|\overline z|
-2|W(G)|
+\sum_{\overline z\in\mathcal U(G)\setminus\mathcal W(G)}(|\overline z|-1)\\
&=&|G|-2|W(G)|-(|\mathcal U(G)|-|\mathcal W(G)|).
\end{eqnarray*}
Since $|\mathcal W(G)|=3|W(G)|$, our desired result follows.
$\qed$

We use $\Psi$ to denote the set of noncyclic groups $G$ satisfying that
there exists an odd prime $p$ such that the following three conditions hold.

(C1) The prime divisors of $|G|$ are $2$ and $p$.

(C2) The subgroup of order $p$ is unique.

(C3) There is no element of order $4$ in $G$.

(C4) Each involution of $G$ is contained in a cyclic subgroup of order $2p$.

\begin{example}
  Let $Z_n$ denote the cyclic group of order $n$. If $m\geq 2$,$n\geq 1$ and $p$ is an odd prime, then
  $$
  \overbrace{Z_2\times\cdots\times Z_2}^m\times Z_{p^n}\in\Psi.
  $$
\end{example}

\begin{prop}
Suppose $G\in\Psi$. Then $|G|=2^mp^n$ for some positive integers $m,n$ and
odd prime $p$. Moreover, the Sylow $2$-subgroup is an elementary abelian
$2$-group and the Sylow $p$-subgroup is a cyclic group.
\end{prop}
\proof The condition (C1) implies that $|G|=2^mp^n$ for some positive integers $m,n$ and odd prime $p$. By (C3),  the Sylow $2$-subgroup is an elementary abelian
$2$-group. It follows from Lemma~\ref{p} and (C2) that the Sylow $p$-subgroup is cyclic.
$\qed$

\begin{lemma}\label{e5}
Let $G$ be a noncyclic group. Then $G\in\Psi$ if
and only if there is a nonidentity element $x$ of $G$ such that
the following conditions hold.

{\rm(i)} All elements in $R\{e,x\}\setminus\{e,x\}$ are involutions.

{\rm(ii)} There exist $r-3$ involutions
in $R\{e,x\}\setminus\{e,x\}$ which are not maximal involutions of $G$,
where $r=\max\{|R\{e,x\}|,4\}$.
\end{lemma}
\proof Suppose $G\in\Psi$. Pick an element $x\in G$ with $o(x)=p$,
where $p$ is an odd prime and $p$ divides $|G|$.
For any element $y\in G$ with $o(y)\geq 3$,    by (C1) and (C3), the prime $p$ divides $o(y)$,
which implies that $\langle x\rangle\subseteq\langle y\rangle$ by (C2),
and so $y\not\in R\{e,x\}$. Hence (i) holds. The condition (C4)
implies that (ii) holds.

Conversely, suppose that there is a nonidentity
element $x$ of $G$ such that (i) and (ii) hold.
Write
$$
R_0=R\{e,x\}\setminus\{e,x\},\qquad
R_1=\{z\mid z\in R_0, z\textup{ is not a maximal involution of }G\}.
$$

We claim that, for any $z\in R_1$,
we have $z\not\in\langle x\rangle$ and there exists
an element $z'\in G$ such that
$\langle z\rangle\cup\langle x\rangle\subseteq \langle z'\rangle$.
In fact, for any $z\in R_1$, since $R_1\subseteq R_0\subseteq R\{e,x\}$,
we have $z\not\in\langle x\rangle$.
By (i) there exists
an element $z'\in G\setminus R_0$
such that $\langle z\rangle\subsetneqq\langle z'\rangle$.
Since $\langle z'\rangle\nsubseteq\langle x\rangle$,
we have $\langle x\rangle\subseteq\langle z'\rangle$.
Hence, our claim  is valid.

By (ii) we get $R_1\neq\emptyset$.
Pick $z_0\in R_1$. By (i) we have $o(z_0)=2$.
By the claim, we have $z_0\not\in\langle x\rangle$ and there is an element $z'_0$ such that
$\langle z_0\rangle\cup\langle x\rangle\subseteq \langle z'_0\rangle$,
which implies that $o(x)$ is odd, $o(z'_0)$ is even and $o(x)$ divides
$o(z'_0)$.
If $o(x)$ is not a prime, there is an even number $m$ with $2<m<o(z'_0)$
such that $m$ divides $o(z'_0)$ and $o(x)$ does not divide $m$,
which implies that any element of order $m$ in $\langle z'_0\rangle$
is in $R_0$, contrary to (i).
Hence $o(x)$ is an odd prime.

Write $p=o(x)$. Then $p$ is an odd prime.
Hence, for any $x'\in G\setminus R\{e,x\}$, we get $x\in\langle x'\rangle$.
Therefore,
the condition (C1), (C2) and (C3) hold.
Note that $R_0$ consists of all involutions in $G$.
In order to prove (C4), we only need to prove $R_0=R_1$.

If $|R\{e,x\}|\leq 3$, then $|R\{e,x\}|=3$ and $|R_0|=1$, which
implies that $R_0=R_1$ by (ii). Now suppose $|R\{e,x\}|\geq 4$.
By (ii), we have
\begin{equation}\label{e4}
0\leq|R_0|-|R_1|\leq 1.
\end{equation}

Write
$m_i=|\{g\mid g\in G,o(g)=p^i\}|$ and $n_i=|\{g\mid g\in G,o(g)=2p^i\}|.$
Let $s$ and $t$ be the maximum
numbers of $i$ such that $m_i\neq 0$ and $n_i\neq 0$, respectively.
By (C1) and (C3), we have
\begin{equation}\label{e3}
\sum_{i=0}^s m_i+\sum_{i=0}^t n_j=|G|.
\end{equation}
Since  $\varphi(p^i)$ divides $m_i$  and $\varphi(2p^i)$ divides $n_i$,
the prime $p$ divides $m_i$  and  $n_i$ for $i\geq 2$, which implies
that $p$ divides $m_0+m_1+n_0+n_1$ by (\ref{e3}).
It is clear that $m_0=1$, $n_0=|R_0|$ and $n_1=|R_1|(p-1)$.
By (C2), we have $m_1=p-1$. So $p$ divides $|R_0|+|R_1|(p-1)$.
It follows from (\ref{e4}) that $R_0=R_1$.
$\qed$

\begin{lemma}\label{e8}
For any $G\in\Psi$ , we have
 $\dim(\mathcal P_G)\geq|G|-|\mathcal U(G)|+1$.
\end{lemma}
\proof By Lemma~\ref{e5}, there exists a nonidentity element $x\in G$ such that $R\{e,x\}\setminus\{e,x\}$ is a collection of involutions. Write $R_0=R\{e,x\}\setminus\{e,x\}$. For each $w\in R_0$, since there is no element of order $4$, by Proposition~\ref{closed} we get $\overline w=[w]=\{w\}$. Let
$$
\mathcal U=\mathcal U(G),\quad\mathcal U_1=\{\overline w\mid w\in R_0\}\cup\{\overline e,\overline x\}\quad\textup{and}\quad A=\bigcup_{\overline z\in\mathcal U_1}\{\overline z\}.
$$

Suppose $S$ is a resolving set of $\mathcal P_G$ with size $\dim(\mathcal P_G)$. By Lemma~\ref{e7}, one gets
\begin{equation}\label{5}
|S\cap A|\geq |\overline x|-1+|\overline e|-1+1=|\overline x|+|\overline e|-1=|A|-|R_0|-1.
\end{equation}
Since $|\mathcal U_1|=|R_0|+2$, by Lemma~\ref{e7} and (\ref{5}), we have
\begin{eqnarray*}
\dim(\mathcal P_G)=|S|=|S\cap(\bigcup_{\overline z\in\mathcal U}\overline z)|&=&|S\cap A|+\sum_{\overline z\in\mathcal U\setminus\mathcal U_1}|S\cap\overline z|\\
&\geq&|A|-|R_0|-1+\sum_{\overline z\in\mathcal U\setminus\mathcal U_1}(|\overline z|-1)\\
&=&|G|-|\mathcal U|+1,
\end{eqnarray*}
as desired.
$\qed$

Now we give a closed formula for the metric dimension of the power graph of a group.

\begin{thm}\label{f1}
Let $G$ be a group.

{\rm(i)} If $G\in\Psi$, then $\dim(\mathcal P_G)=|G|-|\mathcal U(G)|+1$.

{\rm(ii)} If $G\not\in\Psi$, then
 $\dim(\mathcal P_G)=|G|-|\mathcal U(G)|+|W(G)|$.
\end{thm}
\proof Write $\mathcal U(G)=\{\overline{x_1},\ldots,\overline{x_t}\}$ and
$$
X=G\setminus\{x_1,\ldots,x_t\},
$$
where $t=|\mathcal U(G)|$.  For any two distinct elements $u_1$ and $u_2$ in $G$, write
$$
R_0\{u_1,u_2\}=R\{u_1,u_2\}\setminus\{u_1,u_2\}.
$$

{\em Claim 1.  If there is an element  of order at least three
in $R_0\{u_1,u_2\}$, then $X \cap R_0\{u_1,u_2\}\neq\emptyset$.}
In fact, if $z_0\in R_0\{u_1,u_2\}$ and $o(z_0)\geq 3$, by Observations~\ref{equi} and \ref{rxy}, we have $[z_0]\subseteq\overline{z_0}\subseteq R_0\{u_1,u_2\}$, and so $[z_0]\setminus\{x_1,\ldots,x_t\}\subseteq X\cap R_0\{u_1,u_2\}$. Since $|[z_0]|=\varphi(o(z_0))\geq 2$, one gets $[z_0]\setminus\{x_1,\ldots,x_t\}\neq\emptyset$. Consequently,  Claim 1 is valid.

{\em Claim 2. If  all elements in $R_0\{u_1,u_2\}$ are involutions, then  $\langle u_1\rangle\subsetneqq \langle u_2\rangle$ or $\langle u_2\rangle\subsetneqq \langle u_1\rangle$.} If $\langle u_1\rangle=\langle u_2\rangle$, then $\overline{u_1}=\overline{u_2}$, and so $R_0\{u_1,u_2\}=\emptyset$, a contradiction. Suppose that $u_1$ and $u_2$ are not adjacent in $\mathcal P_G$. On one hand, for any $z\in R_0\{u_1,u_2\}$, we conclude that $z$ is adjacent to one of $u_1$ and $u_2$ and not adjacent to the other. Without loss of generality, assume that $z$ is adjacent to  $u_1$. Since $o(z)=2$ and $u_1\neq e$, we have $z\in \langle u_1\rangle$, and so $o(u_1)\geq 4$. On the other hand, since any element of $[u_1]$ is not adjacent to $u_2$ in $\mathcal P_G$,  we have $[u_1]\setminus\{u_1\}\subseteq R_0\{u_1,u_2\}$, which implies that
$o(u_1)=2$, a contradiction. Hence, Claim 2 is valid.

 (i)  By Lemma~\ref{e5}, there exists a nonidentity element $x\in G$ such that $R_0\{e,x\}$ is a collection of involutions.  Pick an element $y_0\in R_0\{e,x\}$. Let
$$
Y=X\cup\{y_0\}.
$$
By Proposition~\ref{closed} and (C3), we have $\overline{y_0}=\{y_0\}$. Then $|Y|=|X|+1=|G|-|\mathcal U(G)|+1$. By Lemma~\ref{e8}, we only need to show that $Y$ is a resolving set of $\mathcal P_G$. Pick any two distinct vertices $u_1$ and $u_2$ in $G\setminus Y$. It suffices to show that
\begin{equation}\label{10}
  Y\cap R_0\{u_1,u_2\}\neq\emptyset.
\end{equation}
If there exists an element  of order at least three in $R_0\{u_1,u_2\}$,  by Claim 1, we have $X \cap R_0\{u_1,u_2\}\neq\emptyset$, which implies that (\ref{10}) holds. Note that $e\not\in R_0\{u_1,u_2\}$. Now suppose that all elements in $R_0\{u_1,u_2\}$ are involutions. By Claim 2, without loss of generality, assume that $\langle u_1\rangle\subsetneqq \langle u_2\rangle$.

In order to prove (\ref{10}), we only need to show that $y_0\in R_0\{u_1,u_2\}$. Suppose for the contrary that $y_0\not\in R_0\{u_1,u_2\}$.

Since $\{u_1,u_2\}\subseteq G\setminus Y\subseteq\{x_1,\ldots,x_t\}$, we have $\overline{u_1}\neq\overline{u_2}$, which implies that $R_0\{u_1,u_2\}\neq\emptyset$. Pick $u_0\in R_0\{u_1,u_2\}$. Then $o(u_0)=2$ and $u_0$ is adjacent to one of $u_1$ and $u_2$ and not adjacent to the other in $\mathcal P_G$. If $u_0$ is adjacent to $u_2$ and not adjacent to $u_1$ in $\mathcal P_G$, then $\langle u_0\rangle\subsetneqq\langle u_2\rangle$ and $u_1\neq e$. Since $u_0$ is the unique involution in the subgroup $\langle u_2\rangle$, we have $R\{u_1,u_2\}\cap\langle u_2\rangle=\{u_1,u_2,u_0\}$, which implies that $u_0$ is a resolving involution of $\langle u_2\rangle$. Let $p$ be an odd prime that divides $|G|$. By Proposition~\ref{ril2},  we have $o(u_1)=p$ and $o(u_2)=2p^m$ for some positive integer $m$. The fact that $o(y_0)=2$ implies that there exists an element $u_3$ of order $2p$ such that $y_0\in\langle u_3\rangle$ by (C4). By (C2), one has $u_1\in\langle u_3\rangle$. Since $y_0\neq u_0$, we have $y_0\not\in\langle u_2\rangle$, and so $u_3\not\in\langle u_2\rangle$.  Therefore, we get $u_3\in R_0\{u_1,u_2\}$, a contradiction. Hence $u_0$ is adjacent to $u_1$ and  not adjacent to $u_2$ in $\mathcal P_G$, which implies that $u_1=e$ and $u_0\not\in\langle u_2\rangle$. By (C4), there exists an element $u_4$ of order $2p$ such that $u_0\in\langle u_4\rangle$. Then $u_4\not\in\langle u_2\rangle.$ Since $u_4\not\in R\{u_1,u_2\}$, we have $u_2\in\langle u_4\rangle$, which implies that $o(u_2)=p$, and so $y_0\in R_0\{u_1,u_2\}$, a contradiction.

(ii)  Write
$$
S=X\cup W(G).
$$
 Observation~\ref{ri} implies that $\overline w=\{w\}$ for each $w\in W(G)$, and so $|S|=|G|-|\mathcal U(G)|+|W(G)|$. By Lemma~\ref{low}, we only need to show that
$S$ is a resolving set of $\mathcal P_G$. Pick any two distinct elements $u_1$ and $u_2$ in $G\setminus S$. It suffices to show that
\begin{equation}\label{9}
S\cap R_0\{u_1,u_2\}\neq\emptyset.
\end{equation}
If there exists an element  of order at least three in $R_0\{u_1,u_2\}$,  by Claim 1, we have $X \cap R_0\{u_1,u_2\}\neq\emptyset$, which implies that (\ref{9}) holds. Note that $e\not\in R_0\{u_1,u_2\}$. Now suppose that all elements in $R_0\{u_1,u_2\}$ are involutions. By Claim 2, without loss of generality, assume that $\langle u_1\rangle\subsetneqq \langle u_2\rangle$.

If $|R_0\{u_1,u_2\}|=1$, then $R_0\{u_1,u_2\}\subseteq W(G)\subseteq S$, and so (\ref{9}) holds. Suppose $|R_0\{u_1,u_2\}|\geq 2$. Since $\langle u_2\rangle$ contains at most one involution, there exists an involution $z_1\in R_0\{u_1,u_2\}\setminus\langle u_2\rangle$. Note that $z_1$ and $u_2$ are not adjacent in $\mathcal P_G$. Then $z_1$ and $u_1$ are adjacent in $\mathcal P_G$, which implies that $u_1=e$ by $z_1\not\in\langle u_1\rangle$. Since $|R\{e,u_2\}|=|R_0\{u_1,u_2\}|+2\geq 4$ and $G\not\in\Psi$, by Lemma~\ref{e5}, there exist two distinct maximal involutions $v_1$ and $v_2$ of $G$ in $R_0\{u_1,u_2\}$. By Lemma~\ref{open}, we have $\overline{v_1}=\overline{v_2}$, and so $\{v_1,v_2\}\cap S\neq\emptyset$, which implies that (\ref{9}) holds.
$\qed$

As a corollary, we compute the metric dimension of the power graph of a cyclic group.

\begin{cor}
Suppose $n=p_1^{r_1}\cdots p_t^{r_t}$, where
$p_1,\ldots,p_t$ are primes  with $p_1<\cdots<p_t$,
and $r_1,\ldots,r_t$ are positive integers.
Let $Z_n$ denote the cyclic group of order $n$.
Then
$$
\dim(\mathcal{P}_{Z_n})=
\left\{
\begin{array}{ll}
n-1, &\textup{if }t=1,\\
n-2r_2, &\textup{if }(t,p_1,r_1)=(2,2,1), \\
n-2r_1, &\textup{if }(t,p_1,r_2)=(2,2,1),\\
n+1-\prod_{i=1}^t(r_i+1), &\textup{otherwise. }
\end{array}\right.
$$
\end{cor}
\proof
If $t=1$, then $n$ is a prime power, which implies
that $\dim(\mathcal P_{Z_n})=n-1$ by Lemma~\ref{complete}.
Now suppose $t\geq 2$. By Propositions~\ref{cycequi}, we have
$$
|\mathcal U(Z_n)|=|\mathcal C'(Z_n)|-1=|\mathcal C(Z_n)|-1=\prod_{i=1}^t(r_i+1)-1.
$$
By Proposition~\ref{ril2}, we have
$$
|W(Z_n)|=
\left\{
\begin{array}{ll}
1, &\textup{if }(t,p_1,r_1,r_2)=(2,2,1,r_2)\textup{ or }(2,2,r_1,1),\\
0, &\textup{otherwise. }
\end{array}\right.
$$
Consequently, Theorem~\ref{f1} (ii) implies that our desired result
follows.
$\qed$

\section*{Acknowledgement}

This research is supported
by National Natural Science Foundation of China (11271047,  11371204).

\end{CJK*}


\begin{thebibliography}{99}

\bibitem{aba} J. Abawajy, A. Kelarev and M. Chowdhury, Power graphs: A survey,
Electron. J. Graph Theory Appl.  1  (2013), 125--147.

\bibitem{RP} R.F. Bailey and P.J. Cameron, Base size, metric dimension and other
invariants of groups and graphs, Bull. London Math. Soc. 43 (2011),
209--242.

\bibitem{bat} C. Bates, D. Bundy, S. Perkins and P. Rowley, Commuting involution graphs for symmetric
groups, J. Algebra, 266 (2003), 133--153.

\bibitem{bol} B. Bollob\'as, Mordern Graph Theory, Springer, New York, 1998.

%\bibitem{bon} F. Bonomo, S. Mattia and G. Oriolo, Bounded coloring of co-comparability graphs and the pickup and delivery tour combination problem, Theoret. Comput. Sci.  412  (2011), 6261--6268.

\bibitem{bos} J. Bos\'ak, The graphs of semigroups, in: Theory of Graphs and Application, Academic Press, New York, 1964, pp. 119--125.

\bibitem{Ca} J. C\'aceres, C. Hernando, M. Mora, I.M. Pelayo, M.L. Puertas, C. Seara
and D.R. Wood, On the metric dimension of Cartesian products of
graphs,  SIAM J. Discrete Math. 21 (2007), 423--441.

\bibitem{cam} P.J. Cameron, The power graph of a finite group II, J.  Group Theory 13 (2010), 779--783.

\bibitem{came} P.J. Cameron, S. Ghosh, The power graph of a finite group, Discrete Math. 311 (2011), 1220--1222.


\bibitem{cha} I. Chakrabarty, S. Ghosh and M.K. Sen, Undirected power graphs of semigroups, Semigroup Forum  78   (2009), 410--426.

%\bibitem{cor} S. Cornelsen and  G. Di Stefano, Track assignment, J. Discrete Algorithms  5  (2007),   250--261.

\bibitem{di} G. Di Stefano, Distance-hereditary comparability graphs, Discrete Appl. Math.  160  (2012), 2669--2680.

\bibitem{doo}   A. Doostabadi, A. Erfanian and A. Jafarzadeh, Some results on the power
graph of groups, The Extended Abstracts of the 44th Annual Iranian Mathe-
matics Conference 27--30 August 2013, Ferdowsi University of Mashhad, Iran,
\url{http://profdoc.um.ac.ir/articles/a/1036567.pdf}

\bibitem{gal} T. Gallai, Transitiv orientierbare graphen, Acta Math. Acad. Sci. Hungar  18   (1967), 25--66.

\bibitem{Ga}  M.R. Garey and D.S. Johnson, Computers and Intractability: A Guide to the Theory of NP-Completeness, Freeman, New York, 1979.

\bibitem{gho} A. Ghouila-Houri, Caract\'erization des graphes non orient\'es dont on peut orienter les ar\v{e}tes de mani\`ere \`a obtenir le graphe d'une relation d'ordre,
 C. R. Acad. Sci. Paris 254 (1962), 1370--1371.

\bibitem{gil} P.C. Gilmore and A.J. Hoffman, A characterization of comparability graphs and of interval graphs, Canad. J. Math. 16 (1964), 539--548.

\bibitem{god} C. Godsil and G. Royle, Algebraic Graph Theory, Springer, New York, 2001.

\bibitem{god2} C. Godsil and G.F. Royle, Cores of geometric graphs, Ann. Comb. 15 (2011), 267--276.

%\bibitem{gol} M.C. Golumbic, Algorithmic Graph Theory and Perfect Graphs, Computer Science and Applied Mathematics, Academic Press, 1980.

\bibitem{gor}    D. Gorenstein, Finite Groups,  Second Edition, Chelsea Publishing Co., New York, 1980.

\bibitem{Ha} F. Harary and R.A. Melter, On the metric dimension of a graph,
 Ars Combin. 2 (1976), 191--195; 4 (1977), 318.

\bibitem{He} C. Hernando, M. Mora, I. M. Pelayo, C. Seara and D. R. Wood, Extremal graph theory for metric dimension and diameter,
Electron. Notes in Discrete Math. 29 (2007), 339--343.

%\bibitem{hib} T. Hibi, The comparability graph of a modular lattice, Combinatorica  18  (1998), 541--548.

\bibitem{liy} N. Iiyori and H. Yamaki, Prime graph components of the simple groups of Lie type over the field of even characteristic, Proc. Japan Acad. Ser. A Math. Sci.  67 (3) (1991), 82--83.

\bibitem{ille} P. Ille and J. Rampon,
Reconstruction of posets with the same comparability graph.
J. Combin. Theory Ser. B  74  (1998), 368--377.

\bibitem{kel1} A.V. Kelarev and S.J. Quinn, A combinatorial property and power graphs of groups,  Contributions to general algebra, 12 (Vienna, 1999),  229--235, Heyn, Klagenfurt, 2000.

\bibitem{kel2} A.V. Kelarev and S.J. Quinn, Directed graph and combinatorial properties of semigroups, J. Algebra 251 (2002), 16--26.

\bibitem{Kh} S. Khuller, B. Raghavachari and A. Rosenfeld, Landmarks in graphs, Discrete Appl. Math. 70 (1996), 217--229.



\bibitem{mir} M. Mirzargar, A.R. Ashrafi and M.J. Nadjafi-Arani, On the power graph of a finite group. Filomat  26  (2012),  1201--1208.

\bibitem{mog} A.R. Moghaddamfar S. Rahbariyan and W.J. Shi, Certain properties of the power graph associated with a finite group, \url{http://arxiv.org/pdf/1310.2032v1.pdf}

%\bibitem{pnu} A. Pnueli, A. Lempel and  S. Even, Transitive orientation of graphs and identification of permutation graphs. Canad. J. Math.  23  (1971), 160¨C175.


\bibitem{Sa} G. Sabidussi, Graph derivates, Math. Z. 76 (1961),  385--401.



\bibitem{Sl} P.J. Slater, Leaves of trees,  Congr. Numer.
14 (1975), 549--559.

%\bibitem{spe} E. Sperner, Ein Satz ber Untermengen einer endlichen Menge, Math. Z. 27 (1928), 554--548.

\bibitem{tam}  T. Tamizh Chelvam and M. Sattanathan, Power graph of finite abelian
groups,  Algebra Discrete Math.  16  (2013),  33--41.

\bibitem{zel} B. Zelinka, Intersection graphs of finite abelian groups, Czechoslovak Math. J. 25  (1975), 171--174.

\end{thebibliography}
\end{document}